\documentclass[12pt]{article}

\usepackage{graphicx}
\usepackage{amsfonts,amsmath,amssymb,bbm}
\usepackage{setspace} 
\def\proof{\noindent{\bf Proof:}\hskip10pt}        
\def\QED{\hfill $\Box$}

\font\tenmath=msbm10 scaled 1200
\font\sevenmath=msbm7 scaled 1200
\font\Fivemath=msbm5 scaled 1200
\newfam\mathfam \textfont\mathfam=\tenmath
\scriptfont\mathfam=\sevenmath \scriptscriptfont\mathfam=\Fivemath

\def \\ { \cr }
\def\R{\mathbb{R}}
\def \1{1 \mkern -6mu 1} 
\def\N{\mathbb{N}}
\def\E{\mathbb{E}}
\def\P{\mathbb{P}}

\def\R{\mathbb{R}}

\def \e{{\rm e}}
\def \d{{\rm d}}
\def \Yp{Y^{(p)}}

\newtheorem{theorem}{Theorem}

\newtheorem{proposition}{Proposition}
\newtheorem{lemma}{Lemma}
\newtheorem{corollary}{Corollary}
\begin{document}

\title{On the non-Gaussian fluctuations\\
 of the giant cluster  for percolation \\ 
 on random recursive trees }

\author{Jean Bertoin\thanks{Institut f\"ur Mathematik, 
Universit\"at Z\"urich, 
Winterthurerstrasse 190, 
CH-8057 Z\"urich}
\thanks{
{email: jean.bertoin@math.uzh.ch}} 
 }

\date{}
\maketitle

\begin{abstract}  
We consider a Bernoulli bond percolation on a random recursive tree of size $n\gg 1$, with supercritical parameter $p_n=1-c/\ln n$ for some $c>0$ fixed. It is known that with high probability, there exists then a unique giant cluster of size $G_n\sim \e^{-c}$, and it follows from a recent result of Schweinsberg \cite{Sch} that
 $G_n$ has non-gaussian fluctuations. We provide an  explanation of this  by analyzing  the effect of percolation on different phases of the growth of recursive trees. 
This alternative approach may be useful for studying percolation on other classes of trees, such as for instance regular trees.
\end{abstract}

  
\begin{section}{Introduction and main result}

A famous result due to  Erd\"os and R\'enyi shows that  Bernoulli bond percolation on the complete graph with $n$ vertices and with parameter $c/n$ for $c>1$ fixed, produces with high probability as $n\to \infty$, a unique giant cluster of size $\Gamma_n\sim \theta(c)n$,  where $\theta(c)$ is  the strictly positive solution to the equation $x+\e^{-cx}=1$. It has been known since the work of Stepanov  \cite{St} that the fluctuations of $\Gamma_n$ are normal, in the sense that 
$$\frac{\Gamma_n-\theta(c)n}{\sqrt{n}} \, \Longrightarrow {\mathcal N}(0,\sigma_c^2),$$ 
where as usual  ${\mathcal N}(0,\sigma_c^2)$ denotes a centered Gaussian variable with variance $\sigma_c^2$,  and $\Rightarrow$ means convergence in law. 
See also e.g. Pittel \cite{Pi} and  Barraez {\it et al.} \cite{BBF} for alternative proofs and refinements. 
Since then, several  results  have appeared in the literature, establishing the asymptotic normality of giant components for various random graph models. We refer in particular to Behrisch {\it et al.} \cite{BCOK},
Bollob\`as and Riordan \cite{BR}, and Seierstad \cite{Se}. 

The motivation of the present work stems from the feature that  the giant cluster resulting from supercritical bond percolation on a large random recursive tree has a much different asymptotic behavior. Recall that a tree on an ordered set of vertices, say $[n]=\{1,\ldots, n\}$, is called {\it recursive} if when rooted at $1$, the sequence of vertices along  any branch from the root to a leaf increases. The terminology comes from the fact that such trees can be constructed recursively, incorporating  each vertex  one after the other in the natural order  to built a growing tree. See Drmota \cite{Dr} for background and further references. 

We denote by $T_n$ a recursive tree picked uniformly at random amongst the $(n-1)!$ recursive trees on $[n]$. Equivalently, $T_n$ can be constructed recursively by creating for   $\ell=1,\ldots, n-1$ an edge between the vertices  $\ell+1$ and $u_{\ell}$, where  $u_{\ell}$ has the uniform distribution on $[\ell]$ and $u_1,\ldots, u_{n-1}$ are independent. Given $T_n$, we then perform a Bernoulli bond percolation with parameter
$$p_n=1-\frac{c}{\ln n}$$
where $c>0$ is some fixed parameter.  It is easy to show that this choice of the percolation parameter corresponds precisely to the supercritical regime, in the sense that with high probability for $n\gg 1$, 
the cluster containing the root is giant with size $G_n\sim \e^{-c} n$. At this point, it may be interesting to briefly sketch the proof of this result, referring to \cite{Be2} for details. 

Pick a vertex $u_n$ uniformly at random in $[n]$, and denote  its distance to the root by $h_n$. Then it is well known that $h_n\sim \ln n$, and  since the first moment of $n^{-1}G_n$ coincides with the probability that $u_n$ is connected to the root, one gets 
$$\E(n^{-1}G_n)=\E\left( \left (1-\frac{c}{\ln n}\right )^{h_n} \right) \sim \e^{-c}. $$ Similarly, if $v_n$ denotes a second uniform vertex  chosen independently of the first, then the  easy fact that the height  of the  branch point of $u_n$ and $v_n$ remains stochastically bounded yields the second moment estimate $\E((n^{-1}G_n)^2)\sim \e^{-2t}$, from which the law of large numbers for $G_n$ follows. 

Since it is also well-known that  $h_n$ is asymptotically normal (see Devroye \cite{De}), this might suggest that the same could also hold for $G_n$. However, it follows from a recent result due to Schweinsberg \cite{Sch}  that this is not the case. In order to give a precise statement, recall that a real-valued random variable  $Z$ has the so-called continuous Luria-Delbr\"uck law when its characteristic function is given by
$$\E(\e^{i \theta Z})= \exp\left( -\frac{\pi}{2} |\theta| -i \theta \ln |\theta| \right)\,, \qquad \theta \in \R\,.$$
This distribution arises in limit theorem for sums of positive  i.i.d. variables in the domain of attraction of a completely asymmetric Cauchy process; see e.g. Geluk and de Haan \cite{GdH}, M\"ohle \cite{Moe},  and further references therein. Its role in the context of  large random recursive trees was observed first by Drmota {\it et al.} Ê\cite{DIMR} and Iksanov and M\"{o}hle \cite{IM}, in relation with a random algorithm for the isolation of the root. See also Holmgren \cite{Holm} and the comments (a) and (b) in the forthcoming Section 2.3. 

\begin{theorem} \label{T1} {\rm (Schweinsberg)}
There is the weak convergence
$$ \left(n^{-1}G_n - \e^{-c} \right)\ln n- c \e^{-c} \,  \ln \ln n  \, \Longrightarrow \, -c\e^{-c}\left( Z + \ln c\right)\,,$$
where the variable  $Z$ has the continuous Luria-Delbr\"uck distribution. 
\end{theorem}

More precisely, Theorem 1.7 in \cite{Sch} is stated in terms of the number of blocks in the Bolthausen-Sznitman coalescent, and the remarkable construction of Goldschmidt and Martin \cite{GM} of the latter based on random cuts on a random recursive tree entails the present statement. The proof of Schweinsberg relies on delicate estimates on the rate of decrease of the number of blocks in the Bolthausen-Sznitman coalescent and on bounds on stable processes, and the purpose of this work is to propose an alternative approach which may provide more intuitive explanations for the anomalous fluctuations of $G_n$. We stress that Theorem 1.7 in \cite{Sch}  is a weak limit theorem for processes while Theorem \ref{T1} only states a one-dimensional convergence result; however the present approach immediately extends to finite-dimensional convergence at the price of slightly heavier notations, and establishing tightness would require some further work.

In the sequel, it will be convenient to agree that   the edges are enumerated naturally in the  order induced by the construction, i.e. the $\ell$-th edge refers to the edge linking the vertex $\ell$ to its parent $u_{\ell-1}$.
Roughly speaking, one can then distinguish three phases of the random dynamics. 
 Because in percolation, each edge is removed with probability $c/\ln n$, the first edges which are removed correspond to an early  phase 
when the growing tree has size of order $ \ln n$. During this phase, only a stochastically  bounded number of edges are removed, and it has been shown  in \cite{BUB} that the percolation clusters corresponding to those edges  will eventually have size of order $n/\ln n$ when the construction is completed. Informally, this is the source of the random fluctuations involving the  Luria-Delbr\"uck variable  in Theorem \ref{T1}. 

There is then an intermediate phase when the tree grows from a size of order $\ln n$ to the size  $\lfloor \ln^4 n \rfloor$,  during which about $c \ln^3 n$ edges are removed. Each of the  percolation cluster born during this phase has only size $o(n/\ln n)$ at the end of the process. However, the  cumulative effect of these clusters  is nonetheless visible and yields the deterministic correction involving the iterated logarithm factor in Theorem \ref{T1}. 

In the final phase when the recursive tree grows from size $\lfloor \ln^4 n \rfloor$ to size $n$,  the root cluster
 grows essentially  regularly, i.e. without inducing further fluctuations. We point out  that the threshold $\ln^4 n$ appearing in this work is somewhat arbitrary, and $\ln^{\alpha}n$ with $\alpha$ close to $4$ would work just as well. It is however  crucial to choose a threshold which is both sufficiently high so that fluctuations are already visible and spread afterwards quite regularly, and  also sufficiently low so that one can estimate the germ of fluctuations with the desired accuracy.

The rest of this paper is mainly devoted to the proof of Theorem \ref {T1}. The starting point of our analysis is that it is useful to incorporate  percolation during the recursive construction of $T_n$, rather than first  constructing completely $T_n$ and then performing percolation.  In Section 2.1, we interrupt the construction of the random recursive tree when it attains the size $k=\lfloor \ln^4n\rfloor$ and perform a percolation on $T_k$ with parameter $p_n$. We obtain a precise estimate of the number $\Delta_k$ of  vertices of $T_k$ which are disconnected from the root; this  can be viewed as the germ of the anomalous fluctuations for $\Delta_n$. In Section 2.2, we resume the construction of the random recursive tree from the size $k=\lfloor \ln^4n\rfloor$ to the size $n$. Using the basic connexion between random recursive trees and Yule processes, we show that the germ of the anomalous fluctuations $\Delta_k$ spread regularly. Section 2.3 contains some miscellaneous comments.
Finally, in Section 3, we briefly point out that the present approach also applies to study the fluctuations of the size of the giant cluster for percolation on regular trees.  
\end{section}

\begin{section} {Proof of Theorem \ref{T1}}
\subsection{The germ of anomalous fluctuations}

Imagine that we interrupt the construction of the random recursive tree when it reaches size $k=\lfloor \ln^4n\rfloor$; plainly this yields 
a random recursive tree of size $k$ which we denote by $T_k$. Our purpose in this section is to estimate precisely the number  of vertices which are disconnected from the root when one performs  a bond percolation on $T_k$ with parameter $p_n$. It is convenient to work with the parameter $k$ rather than $n$, that is  we write $p_n=q_k$ and note 
the change in the asymptotic regime of the percolation parameter as a function of the size of the tree :
\begin{equation}\label{E7}
q_k= 1-c k^{-1/4} + o(k^{-1}).
\end{equation}
 We shall  establish the following limit theorem in law.

\begin{proposition}\label{P1}  As $k\to \infty$, there is the weak convergence
$$k^{-3/4} \Delta_k - \frac{3}{4}c  \ln k  \, \Longrightarrow \,c\left( Z+ \ln c\right) $$
where $Z$ has the  continuous Luria-Delbr\"uck distribution
\end{proposition} 

The rest of this section is devoted to the proof of Proposition \ref{P1}. Our guiding line is similar to that  in \cite{Be1}, although the percolation parameter there  had a different asymptotic behavior. Namely we shall work with  a continuous-time version of percolation in which edges are removed independently one of the others at a given rate, and consider the process  that counts the number of vertices which are disconnected from the root as time passed. It suffices to focus on the cuts made to the root-cluster, and we interpret the latter as a continuous-time version of a random algorithm introduced by Meir and Moon \cite{MM1, MM2}  for the isolation of the root. In turn, this enables us to use a coupling due to Iksanov and M\"{o}hle \cite{IM} and reduces the problem to the analysis of the asymptotic behavior of a remarkable random walk in the domain of attraction of a completely asymmetric Cauchy distribution. 

We shall follow the route sketched  above, but in the reverse order for an easier articulation of  the argument. 
To start with, we recall briefly an asymptotic result  on a random walk which plays an important role in the study of the isolation of the roof for random recursive trees.
Let $\xi$ denote an integer-valued random variable with law
\begin{equation}\label{EQL}
 \P(\xi=j)=\frac{1}{j(j+1)}\,, \qquad j=1,2, \ldots
\end{equation}
We consider the random walk 
$$S_\ell=\xi_1+\cdots + \xi_\ell\,, \qquad \ell \in \N\,,$$ where the $\xi_i$ are independent copies of $\xi$. According for instance to \cite{GdH}, there is the weak convergence
\begin{equation}\label{E1}
\ell^{-1} S_{\ell} - \ln \ell   \, \Longrightarrow \, Z
\end{equation}
where $Z$ has  the continuous Luria-Delbr\"{u}ck law. 

Iksanov and M\"ohle have pointed at a useful coupling  which connects the preceding random walk and an algorithm of Meir and Moon \cite{MM1, MM2} for the isolation of the root. 
Following Meir and Moon, we define recursively  a decreasing sequence of subtrees 
$$T_k=T_k(0)\supset T_k(1)\supset \ldots$$ 
as follows. We pick a first  edge uniformly at random amongst the $k-1$ edges of $T_k$ and remove it, which disconnects $T_k$ into two subtrees. We denote by $T_k(1)$ the subtree which contains the root and imagine that the subtree which does not contain the root is set aside. Then we pick second edge uniformly at random in $T_k(1)$, remove it. We write $T_k(2)$ for the subtree which contains the root, set the other subtree aside, and iterate in an obvious way until the root is finally isolated. 
We write 
$$D_k(\ell)=k-|T_k(\ell)|$$ 
for the number of vertices which have been disconnected from the root after $\ell$ steps, that is the sum of the sizes of the subtrees which have been set aside. 

 Iksanov and M\"{o}hle \cite{IM} have observed that one may couple the random walk $S$ and the algorithm of isolation of the root described above in such a way that there is the identity
\begin{equation}\label{E0}
D_k(\ell) = S_{\ell}Ê\qquad \hbox{for all } \ell< N(k)\,,
\end{equation}
where $N(k)=\min\{\ell: S_{\ell}\geq k\}$ denotes the first passage time of the random walk above level $k$. 
See also Lemma 2 in \cite{Be1} for a statement tailored for our needs. It follows easily that, provided that the number of removed edges is relatively small, $D_k(\ell)$ grows  nearly linearly. Here is a rather crude bound which will be however sufficient for our purpose.

\begin{lemma}\label{L1} Suppose that $\ell= \ell(k)$ fulfills  $1\ll  \ell \leq k\ln^{-2}k$. 
Then 
$$\lim_{k\to \infty} \frac{D_k(\ell)}{\ell \ln^2 \ell}= 0 \qquad \hbox{in probability}\,.$$

\end{lemma}

\proof Indeed, it follows immediately from \eqref{E1} that
\begin{equation}\label{E2}
\lim_{\ell\to \infty} \frac{S_{\ell}}{\ell \ln^2\ell}=0\qquad \hbox{in probability.}
\end{equation}
In particular the assumption $\ell \leq k\ln^{-2}k$ ensures that $\ell< N(k)$ with high probability, so we can use the coupling \eqref{E0} of Iksanov and M\"{o}hle. Then \eqref{E2} is precisely our statement. \QED

We now turn our attention  to a continuous time version of bond percolation on $T_k$. We equip each of its $k-1$ edges with an independent exponential variable with parameter $k^{-1/4}$, and remove each edge at the time given by this variable.  We define
$$t_k= - k^{1/4}\ln q_k\,,$$
so that the probability that a given edge has not yet been removed at time $t_k$ is $\exp(- k^{-1/4} t_k)=q_k$, and 
 the configuration observed at time $t_k$ is thus precisely that resulting from a bond percolation on $T_k$ with parameter $q_k$. 
Note also from \eqref{E7} that
\begin{equation} \label{E3}
t_k= c + O(k^{-1/4})\,.
\end{equation}

As we are only  interested in the number of vertices which have been disconnected from the root at time $t_n$, we may focus on the evolution of the cluster which contains the root. 
Plainly, if we write $\rho_k(\ell)$ for the instant when the $\ell$-th edge is removed from the root-cluster in this continuous-time percolation, then the root-cluster at time 
$\rho_k(\ell)$ can be identified  as $T_k(\ell)$, the subtree obtained by the isolation of the root algorithm after $\ell$ steps. 
We will need the following bounds.

\begin{lemma}\label{L2} Take any $\alpha \in (1/2, 3/4)$. Then 
$$\lim_{k\to \infty} \P\left( \rho_k\left( c k^{3/4} - k^{\alpha} \right) \leq t_k \leq \rho_k\left( c k^{3/4} + k^{\alpha} \right) \right) =1.$$

\end{lemma}

\proof It should  be clear from the  dynamics of continuous-time percolation and elementary properties of independent exponential variables that
 $\rho_k(\ell)$ can be expressed in the form 
 $$\rho_k(\ell) = \sum_{j=0}^{\ell -1} \frac{k^{1/4}}{k-D_k(j)-1}\varepsilon_j$$
 where $\varepsilon_0, \varepsilon_1, \ldots$ is a sequence of i.i.d. standard exponential variables, which is further independent of the algorithm of isolation of the root (the denominator in the fraction above is the number of edges of $T_k(j)$, and for $j$ exceeding the number of steps needed to isolate the root, the general term of the series becomes infinite by convention). 
 
 We take first $\ell =  c k^{3/4} - k^{\alpha}$ and use the obvious lower-bound 
  $$\rho_k(c k^{3/4} - k^{\alpha}) \geq  k^{-3/4} \sum_{j=0}^{c k^{3/4} - k^{\alpha} -1} \varepsilon_j.$$
  Elementary arguments based on the computation of first moment and variance show that the right-hand side can be bounded from below by $c - 2k^{\alpha-3/4}$  with high probability as $k\to \infty$ (note that $\alpha > 3/8$). 
 
Similarly, we then take $\ell =  c k^{3/4} + k^{\alpha}$ and use Lemma \ref{L1} to see that with high probability for $k\gg 1$, there is  the upper-bound 
  $$\rho_k(c k^{3/4} + k^{\alpha}) \leq \frac {k^{1/4}}{k-(c k^{3/4} + k^{\alpha})\ln^2k} \sum_{j=0}^{c k^{3/4} + k^{\alpha} -1} \varepsilon_j.$$
  Again, the sum in the right hand side is easily bounded from above by $c k^{3/4} + 2k^{\alpha}$ with high probability. 
  On the other hand, the quotient is bounded from above by $k^{-3/4}(1+2c k^{-1/4}\ln^2k)$ whenever $k$ is sufficiently large. Putting the pieces together and recalling that $3/4-\alpha < 1/4$, we get that 
  $$\rho_k(c k^{3/4} + k^{\alpha})\leq c+3k^{\alpha-3/4}$$ with high probability for $k\gg 1$. We conclude the proof with an appeal to \eqref{E3}. \QED 

We are now able to establish Proposition \ref{P1}. 

\proof Lemma \ref{L2} and an argument of monotonicity show that  for any $\alpha \in (1/2, 3/4)$,  the bounds
$$D_k( c k^{3/4} - k^{\alpha} ) \leq \Delta_k \leq D_k( c k^{3/4} + k^{\alpha} )$$
hold with high probability.
On the other hand, we see from \eqref{E2} that 
$$\lim_{k\to \infty} k^{-3/4} S_{ k^{\alpha} }  =0 \qquad \hbox{in probability,}$$
and we also know from \eqref{E1} that
$$k^{-3/4} S_{ c k^{3/4} } - \frac{3}{4}c  \ln k  \ \Longrightarrow\ c\left( Z+ \ln c\right).$$
It follows that 
$$k^{-3/4} S_{ c k^{3/4} \pm k^{\alpha} } - \frac{3}{4}c  \ln k \ \Longrightarrow\ c\left( Z+ \ln c\right),$$
and an appeal to  the coupling \eqref{E0} completes the proof. \QED

\subsection  {The spread of anomalous fluctuations}
We now recall the well-known  connexion between random recursive trees and the Yule process. Imagine that  once the tree $T_{\ell}$ of size $\ell=1, \ldots, n-1$ has been constructed, the vertex $\ell+1$ is incorporated after an exponential time with parameter $\ell$, and then connected by a new edge to some vertex $u_{\ell}\in[\ell]$ which is picked uniformly at random and independently of the exponential waiting time. We further mark that edge with probability $1-p_n=c/\ln n$, independently of the preceding events. A mark on an edge indicates that this edge will be removed when percolation is performed; equivalently it can also be interpreted as a mutation occurring in the population. 
The dynamics described above are those of a Yule process with unit rate of birth per individual and with rare neutral mutations  which affect each birth event with probability $c/\ln n$, independently of the other  birth events.

In this section, we begin our observation of this process with rare neutral mutations once it has reached the size $k=\lfloor \ln^4n\rfloor$.   We thus write 
$Y=(Y_t)_{t\geq 0}$ for a standard Yule process started from $Y_0=k$, and consider the time 
$$\tau(n)=\inf\{t\geq 0: Y_t=n\}$$
at which it hits $n$. 
Equivalently, $\tau(n)$ is the time  needed to complete the construction of $T_n$ from $T_k$. 
We shall  first estimate this quantity. 

In the sequel, we shall often use the notation 
$$A_n= B_n+o(f(n))\qquad \hbox{ in probability,} $$where $A_n$ and $B_n$ are two sequences of real random variables and $f: \N\to (0,\infty)$ a function,
to indicate that $\lim_{n\to \infty}|A_n-B_n|/f(n)=0$ in probability. 
\begin{lemma}\label{L3} We have
$$  \e^{\tau(n)}= \frac{n}{\ln^4 n}  + o(1/\ln n) \qquad \hbox{in probability.}$$
\end{lemma}

\proof Elementary properties of Yule processes (see, e.g., Equation (6)  in \cite{dlF}) show that
$$\lim_{n\to \infty} \P\left(\left| n \e^{-\tau(n)}  -k\right | > k^{\alpha}\right) = 0\,.$$
for all $\alpha >1/2$. In particular, for $\alpha < 3/4$, this yields
$$ n\e^{-\tau(n)} =   \left(1+ o(\ln^{-1}n) \right)  \ln^{4}n  \qquad \hbox{in probability.}$$
Our claim follows. 
 \QED

An individual in the population corresponds to a vertex on the tree, and vice-versa. It is called a mutant if its ancestral lineage (i.e. its branch to the root) contains at least one mark, and a clone otherwise. 
The population of size $k$ at the time when we start our observation consists in $\Delta_k$ mutants and $k-\Delta_k$ clones.  We focus on the clone population and write $Y^{\rm (c)}=(Y^{\rm (c)}_t, t\geq 0)$ for the process that counts the number of clones as time passes. Because each clone gives birth to a clone child at rate $p_n$ and independently of the other clones, $Y^{\rm (c)}$ is a Yule process with reproduction rate $p_n$ per individual, and started from 
$Y^{\rm (c)}_0=k-\Delta_k$.

 In this framework, it should be plain that   the size $G_n$ of the root-cluster of $T_n$ after percolation with parameter $p_n$, coincides with the 
number of clone individuals at the time when the total population generated by the Yule process reaches $n$, i.e. there is the identity
\begin{equation}\label{E6}
G_n=Y^{\rm (c)}_{\tau(n)}\,.
\end{equation}
One readily get the following estimate.

\begin{lemma}\label{L4}
 We have
$$  G_n=  \e^{p_n\tau(n)} \left(\ln^4n-\Delta_k\right) + o(n/\ln n) \qquad \hbox{in probability,}$$
  \end{lemma}

\proof  Again from the basic estimate  of Equation (6) in \cite{dlF} and the fact that $Y^{\rm(c)}_0\leq k$, we have for any $\alpha >1/2$ that 
$$ \lim_{n\to \infty} \P\left(\left|  \e^{-p_n\tau(n)} Y^{\rm(c)}_{\tau(n)} - Y^{\rm(c)}_0\right | >   k^{\alpha}\right)= 0\,.$$
We choose $\alpha < 3/4$ and deduce that 
$$Y^{\rm(c)}_{\tau(n)}= \e^{p_n\tau(n)} \left(\ln^4n-\Delta_k\right) + \e^{p_n\tau(n)} o(\ln^3n)\,,  \qquad \hbox{in probability.}$$
Since $p_n\leq 1$, we see from Lemma \ref{L3} that 
$$\e^{p_n\tau(n)} o(\ln^3n)= o(n/\ln n)\,,$$
and our claim follows from \eqref{E6}. \QED

We have now all the ingredients to establish Theorem \ref{T1}

\noindent{\bf Proof of Theorem \ref{T1}: \hskip10pt}      
First, it is convenient to apply Skorokhod's representation theorem and assume that the weak convergence in Proposition \ref{P1} holds in fact almost surely. This enables us to write 
$$\ln^4n- \Delta_k =\ln^4n - \ln^3n\left(3c\ln \ln n + c(Z+\ln c)\right) +o(\ln^3n)\,,  \qquad \hbox{almost surely,}$$
and then to re-express Lemma \ref{L4} in the form 
$$G_n= \e^{p_n\tau(n)}\left(\ln^4n - \ln^3n\left(3c\ln \ln n + c(Z+\ln c)\right)\right)+  o(n/\ln n)\,,  \qquad \hbox{in probability.}$$

We next  note from Lemma \ref{L3} that 
$$\e^{p_n\tau(n)} = \left( \frac{n}{\ln^4 n}  + o(1/\ln n)  \right) ^{1-c/\ln n}\,,  \qquad \hbox{in probability,}$$
and it follows from a couple of lines of calculations that 
 $$\e^{p_n\tau(n)} = \e^{-c} \frac{n}{\ln^4n} + 4c \e^{-c} n \frac{\ln \ln n}{\ln^5n} + o(\ln^{-1}n)\,,  \qquad \hbox{in probability.}$$
 Another line of calculation yields
$$G_n=  \e^{-c}n + c\e^{-c} n \frac{\ln\ln n}{\ln n} - c\e^{-c} \frac{ n}{\ln n}\left( Z+\ln c\right) + o(n/\ln n) \,,  \qquad \hbox{in probability,}$$
which completes the proof. \QED

It may be interesting to point out that the same technique can be applied to estimate the descent of the initial  mutant population. Specifically, the sub-population that stems from the $\Delta_k$ mutants at the initial time is described by a Yule process $Y^{\rm (m)}$ with unit birth rate per individual and  started from $Y^{\rm (m)}_0 = \Delta_k$. 
It should be plain that $Y^{\rm(m)}_{\tau(n)}$ coincides with $\Delta_{k,n}$, the number of vertices $i\in [n]$ such that, on the branch from $i$ to the root  $1$,  at least one edge with label at most $k$  is removed when performing percolation. 
 From the same argument as in Lemma \ref{L4}, one can check that  for any $1/2<\alpha <1$,
 $$\lim_{n\to \infty} \P\left(\left| \Delta_{k,n} -  \e^{\tau(n)} \Delta_k\right | > \e^{\tau(n)}  k^{\alpha}\right) = \lim_{n\to \infty} \P\left(\left|  \e^{-\tau(n)} Y^{\rm(m)}_{\tau(n)} - Y^{\rm(m)}_0\right | >   k^{\alpha}\right)= 0\,.$$
It follows that 
\begin{equation}Ê\label{EM}
\Delta_{k,n} = \e^{\tau(n)} \Delta_k + o(n/\ln n)\qquad \hbox{in probability.}
\end{equation}

On the one hand, recall from Lemma \ref{L3} that $\e^{\tau(n)} k^{\alpha} = o(n /\ln n )$ in probability.
On the other hand, combining Proposition \ref{P1} and Lemma \ref{L3}, we get 
$$\frac{\ln n}{n}Ê\e^{\tau(n)} \Delta_k - 3c \ln\ln n 
\Longrightarrow \,c\left( Z+ \ln c\right),$$
and we conclude that 
\begin{equation} \label{EL}
\frac{\ln n}{n}Ê\Delta_{k,n} - 3c\ln \ln n \, \Longrightarrow \, c\left( Z + \ln c\right)\,.
\end{equation}
More precisely, this weak convergence holds jointly with that in Theorem \ref{T1}, as we can  see from Lemma \ref{L4} and \eqref{EM}

\subsection {Miscellaneous comments}  
For the purpose of this section, 
 it is convenient to rewrite Theorem \ref{T1}
in terms of  $\Delta_n=n-G_n$, the number  of vertices in $T_n$ which are disconnected from the root after performing a bond percolation with parameter $p_n$. 
We have then 
\begin{equation}\label{E7}
 \left(n^{-1}\Delta_n - (1-\e^{-c})\right) \ln n+ c \e^{-c} \,  \ln\ln n  \, \Longrightarrow \, c\e^{-c}\left( Z + \ln c\right)\,.
\end{equation}

We also introduce a standard Luria-Delbr\"uck variable $Z_m$ with parameter $m>0$,  which has generating function
$$\E\left( s^{Z_m}\right) = (1-s)^{m(1-s)/s}\,, \qquad 0\leq s \leq 1\,.$$
Recall that as $m\to \infty$, there is the weak convergence
\begin{equation}\label{*}
\frac{Z_m}{m}-\ln m  \, \Longrightarrow \,  Z
\end{equation}
where $Z$ has the continuous Luria-Delbr\"uck distribution. See Pakes \cite{Pa} or Theorem 4.1 in M\"ohle \cite{Moe}. 
 
{\bf (a)}  It has been argued  that for certain populations models with a small rate of neutral mutation, the number of mutants  has a Luria-Delbr\"uck law ; see Section 2 in Kemp \cite{Ke} and references therein. In this setting, the parameter is given by $m=gN(a+g)$ where $N$ is the total population size, $a$ the rate of birth of clones, and $g$ the rate of birth of new mutants. We stress however that, as pointed out by Kemp, the models leading to these  Luria-Delbr\"uck laws `{\sl involve simplifying assumptions that leave realism somewhat in doubt}'.

In our framework, interpreting the recursive construction of $T_n$ as a Yule process and percolation as rare neutral mutations, this suggests that the number $\Delta_n$ of vertices disconnected from the root  might have a distribution close to the Luria-Delbr\"uck law with parameter $m=cn/\ln n$. If we write $\Delta'_n=Z_m$ for the latter, then \eqref{*} yields the weak convergence 
$$\left(n^{-1}\Delta'_n - c \right)\ln n  + c  \ln\ln n  \, \Longrightarrow \, c\left( Z + \ln c\right)\,.$$
This resembles \eqref{E7}, but  with fundamental discrepancies. Note in particular  that for $c>1$, the estimation above would imply that for $n\gg 1$, the mutant population is close to $cn$, a quantity strictly larger than the total population! It is therefore unlikely that Theorem \ref{T1} could established rigorously from such arguments.

{\bf (b)} If we write  $C_{1,n}\geq  C_{2,n}\geq  \ldots$ for  the sequence of the sizes of the percolation clusters disconnected from the root and ranked in the decreasing order, then 
there is clearly  the identity
\begin{equation}\label{E8}
\Delta_n = \sum_i C_{i,n}\,.
\end{equation}
 Theorem 1 in \cite{Be1} states that for every fixed integer $j$, 
\begin{equation}\label{E9}
\left(\frac{\ln n}{n}C_{1,n}, \ldots, \frac{\ln n}{n}C_{j,n}\right)  \, \Longrightarrow \,  
({\tt x}_1, \ldots, {\tt x}_j)
\end{equation}
where ${\tt x}_1>{\tt x}_2>\ldots$ denotes the sequence of the atoms of a Poisson random measure
on $(0,\infty)$ with intensity $c\e^{-c} x^{-2}\d x\,.$
It is certainly tempting to expect that the finite dimensional convergence \eqref{E9} might be reinforced and then yield \eqref{E7} via \eqref{E8}.

An obvious obstacle is that the series $\sum{\tt x}_i$ diverges a.s.; however this can be circumvented by considering 
$$X_n:= \sum_i  \left \lfloor \frac{n    }{\ln n} \, {\tt x}_i\right \rfloor.$$
The reason for taking integer parts above is of course because cluster sizes are integers. Note that this limits {\it de facto} the sum to atoms  such that ${\tt x}_i\geq n^{-1}\ln n$, and then $X_n < \infty$ a.s. More precisely,  by the elementary mapping theorem for Poisson random measures, $ \lfloor  {\tt x}_1  n \ln^{-1} n\rfloor , \ldots $ can be viewed as the sequence of atoms of a Poisson random measure on $\N$ with intensity $m\mu$ where $m=c\e^{-c} n \ln^{-1} n$ and $\mu$ is the probability measure given by $\mu(k)=k^{-1}-(k+1)^{-1}$. 
Thus $X_n=Z_m$ has the Luria-Delbr\"uck law with parameter $m$, see Section 3 in M\"ohle \cite{Moe}. 

As a consequence of \eqref{*}, there is the weak convergence 
$$ \left(n^{-1}X_n - c\e^{-c} \right) \ln n+ c \e^{-c} \ln\ln n \, \Longrightarrow \, c\e^{-c}\left( Z + \ln c - c\right)\,.
$$
This again resembles  \eqref{E7}, in particular one captures the deterministic  correction involving the iterated logarithm, and the random fluctuations are the same up-to a constant. 
This corroborates the fact that the fluctuations for the size of the giant component are chiefly due to the largest percolation clusters. 
However, this fails to give the correct first order ($c \e^{-c}$ instead of $1-\e^{-c}$),  showing that Theorem \ref{T1} cannot be derived from 
weak limits theorems as \eqref{E9}.

{\bf (c)} We now conclude this section by pointing out that the problem considered in this work could also have been formulated in terms of an urn model {\it \`a la} Polya-Hoppe. 
Indeed, the recursive construction of $T_n$ together with marks on edges corresponding to percolation can also be described as follows. We start with an urn containing a single red ball (the root), and at each step, we add either a red ball or a black ball according to the following random algorithm. With probability $c/\ln n$, we add a black ball to the urn, and with probability $p_n=1-c/\ln n$, we pick a ball uniformly at random in the current contain of the urn, and then replace it into the urn together with a new ball of the same color. Then $\Delta_n$ is the number of black balls when the urn contains exactly $n$ balls, and \eqref{E9} thus gives then a precise limit theorem for the proportion of black balls. 

\end{section}

\begin{section}{Percolation on a regular tree}
The purpose of this section is to point out that the approach used in the proof of Theorem \ref{T1} can be also applied to study percolation on other classes of trees; 
 we shall focus here on the simplest case, namely regular trees. Specifically,  let $d\geq 2$ be a fixed integer, consider the rooted infinite $d$-regular tree (i.e. each vertex has outer-degree $d$) and perform a Bernoulli bond percolation with parameter 
  $$p'_h=\exp(-c/h)\,,$$
 where $c>0$ is fixed and $h\in\N$.
 Observe that the probability that a given vertex at height $h$ has been disconnected from the root equals $1-\e^{-c}$. Since there are $d^h$ vertices at height $h$, 
first and second moments calculations as explained in the Introduction readily show that the number $\nabla_h$ of vertices at height $h$ which have been disconnected from the root fulfills 
 $$\lim_{h\to \infty} d^{-h}\nabla_h = 1-\e^{-c}\qquad \hbox{in probability.}$$

We are interested in the fluctuations of $\nabla_h$. In this direction, it is convenient to use the notation $\log_d x = \ln x / \ln d$ for the logarithm with base $d$ of $x>0$, and $y=\lfloor y\rfloor +\{y\}$ for the decomposition of a real number $y$ as the sum of its integer and fractional parts. 
We introduce for every $b\in[0,1)$ and $x>0$
$$\bar \Lambda_b(x)= \frac{d^{\lfloor b-\log_dx \rfloor+1}}{d-1}.$$
The function $\bar \Lambda_b$ decreases and can be viewed as the tail of a measure $\Lambda_b$
on $(0,\infty)$.  Clearly  $\bar \Lambda_b(x) \asymp x^{-1}$, and it follows in particular  that $ \Lambda_b$ fulfills the integral condition
$$\int_{(0,\infty)}(1\wedge x^2) \Lambda_b(\d x)<\infty\,.$$
This enables us to  introduce a spectrally positive L\'evy process $L_b=(L_b(t))_{t\geq0}$ with Laplace exponent
$$\Psi_b(a)= \int_{(0,\infty)} (\e^{-a x}-1 + ax {\bf 1}_{\{x<1\}}) \Lambda_b(\d x)\,,$$
that is
$$\E\left(\exp(-aL_b(t))\right)= \exp\left\{t \Psi_b(a)\right\}\,,\qquad a\geq 0\,.$$
We stress that a similar process arises in relation with limit theorems for the number of random records on a complete regular tree; see Janson \cite{Jan}. 

We are now able to state the following analog of Theorem \ref{T1} (or rather of  Equation \eqref{E7}).

\begin{theorem} \label{T2} 
In the regime where $h\to \infty$  with  $\{\log_d h\}\to b\in[0,1)$, here is the weak convergence
$$ h\left(d^{-h}\nabla_h - (1-\e^{-c})\right) + c \e^{-c}  \log_d h  \, \Longrightarrow \,
\e^{-c}\left( L_b(c)+ c b \right)\,.$$
\end{theorem}

We now prepare the ground for the proof of Theorem \ref{T2}. The main technical issue is to analyze 
the birth of fluctuations of $\nabla_h$, their propagation being  then an easy matter.

For every $k\in\N$, we enumerate the $d^k$ vertices at  height $k$, and for $i=1, \ldots, d^k$, we write $\eta^{(h)}_{k,i}$ for the total number of edges on the branch from the $i$-th vertex  to the root which have been deleted after percolation with parameter $\e^{-c/h}$. So each $\eta^{(h)}_{k,i}$ has the binomial distribution with  parameter $(k,1-\e^{-c/h})$ and $\eta^{(h)}_{k,i}=0$ if and only if that vertex is still connected to the root. We write
$$\nabla^{(h)}_k= \sum_{i=1}^{d^k} {\bf 1}_{\{\eta^{(h)}_{k,i}\geq 1\}}$$
for the number of vertices at height $k$ which are disconnected from the root, in particular for $k=h$, $\nabla_h=\nabla^{(h)}_h$. We  also
set
$$\Sigma^{(h)}_k = \sum_{i=1}^{d^k} \eta^{(h)}_{k,i}\,.$$
 Clearly, $\nabla^{(h)}_k \leq \Sigma^{(h)}_k$, and the purpose of the next lemma is to point out that these two quantities as close when $k\ll h$. This will be useful in the sequel  as the distribution of $\Sigma^{(h)}_k$ is easier to estimate than that of $\nabla^{(h)}_k$. 

\begin{lemma} \label{L5} We have
$$\E(\nabla^{(h)}_k) = d^k\left(1-\e^{-ck/h}\right) 
\quad \hbox{ and } \quad 
\E(\Sigma^{(h)}_k) = kd^k\left(1-\e^{-c/h}\right).$$
\end{lemma}
\proof The probability that a given vertex at level $k$ has been disconnected from the root equals $1-\P(\eta^{(h)}_{k,i}=0)=1-\e^{-ck/h}$, and as there are $d^k$ vertices at that level, the first assertion is obvious.
Next, consider the edges at height $\ell\in\N$, and for $j=1, \ldots, d^{\ell}$, write $\varepsilon_{\ell,j}=1$ if the $j$-th edge is removed after percolation and $\varepsilon_{\ell,j}=0$ otherwise. So the $\varepsilon_{\ell,j}$ are i.i.d. Bernoulli variables with parameter $1-\e^{-c/h}$ and $\eta^{(h)}_{k,i}=\sum_{\ell=1}^k \varepsilon_{\ell, \ell(i)}$ where $\ell(i)$ denotes the ancestor of $i$ at level $\ell$. Because there are exactly $ d^{k-\ell}$ vertices at level $k$ whose branch to the root passes through a given edge at height $\ell$, there is  the identity
\begin{equation}\label{E14}
\Sigma^{(h)}_k = \sum_{\ell=1}^k \sum_{j=1}^{d^{\ell}} d^{k-\ell} \epsilon_{\ell,j}\,,
\end{equation}
and this  yields our second assertion.  \QED

We next analyze the asymptotic behavior of $\Sigma^{(h)}_k$ in appropriate regimes.

\begin{lemma} \label{L6} In the regime where $h\to \infty$ with $k^2\ll h \ll d^k$ and $\{\log_d h\}\to b\in[0,1)$, there is the weak convergence 
$$h d^{-k} \Sigma^{(h)}_k - c(k-\lfloor \log_d h\rfloor)  \ \Longrightarrow \ L_b(c). $$
\end{lemma}
\proof We compute the Laplace transform of $\Sigma^{(h)}_k $ from the identity \eqref{E14} and get
$$\E\left(\exp(-a \Sigma^{(h)}_k )\right) = \exp\left( - \sum_{\ell=1}^k d^{\ell} \ln \left( \e^{-c/h} + (1-\e^{-c/h})\e^{-a d^{k-\ell}}\right) \right)\,.$$
As a consequence, the cumulant of $h d^{-k} \Sigma^{(h)}_k$, 
$$\kappa^{(h)}_k(a)=-\ln \E\left(\exp(-a h d^{-k} \Sigma^{(h)}_k )\right)$$
can be expressed as
$$\kappa^{(h)}_k(a)= \sum_{\ell=1}^k d^{\ell} \ln \left(1-  (1-\e^{-c/h})(1-\e^{-a hd^{-\ell}})\right).$$
In the regime where $h\to \infty$ with $k^2\ll h$, we get 
\begin{eqnarray*}
\kappa^{(h)}_k(a)&=& c \sum_{\ell=1}^k \frac{d^{\ell}}{h} (1-\e^{-a hd^{-\ell}})+ o(1)\\
&=& c\int_{(0,\infty)} (1-\e^{-a x}- ax {\bf 1}_{\{x\leq 1\}})\Pi^{(h)}_k(\d x)+ ac\int_{(0,1]} x \Pi^{(h)}_k(\d x)+ o(1)\,,
\end{eqnarray*}
where the measure $\Pi^{(h)}_k$ is given by
$$\Pi^{(h)}_k=\sum_{\ell=1}^k \frac{d^{\ell}}{h} \delta_{h d ^{-\ell}}\,.$$

One has immediately
$$\int_{(0,1]} x \Pi^{(h)}_k(\d x)= k-\lfloor \log_dh\rfloor\,.$$
On the other hand, the tail $\bar \Pi^{(h)}_k(x)=\Pi^{(h)}_k((x,\infty))$ is given for $h\ll d^k$ by 
$$\bar \Pi^{(h)}_k(x)= \sum_{\ell=1}^k \frac{d^{\ell}}{h}{\bf 1}_{\{hd^{-\ell}>x\}}=\sum_{\ell=1}^{\lfloor \log_d(h/x)\rfloor} \frac{d^{\ell}}{h}\,.$$
We write the quantity above as
$$h^{-1}\, \frac{d^{\lfloor \log_d(h/x)\rfloor+1}-d}{d-1}= \frac{d^{\lfloor \{\log_d h\}-\log_dx \rfloor+1}-dh^{-1}}{d-1}\,,$$
so that in the regime $h\to \infty$ with $\{\log_d h\}\to b$, we have
$$\lim \bar \Pi^{(h)}_k(x)=  \frac{d^{\lfloor b-\log_dx \rfloor+1}}{d-1}=\bar \Lambda_b(x)\,.$$ 

It is now easy to conclude that in the regime of the statement, 
$$\lim \left( \kappa^{(h)}_k(a) -ac (k-\lfloor \log_dh\rfloor)\right)=-c\Psi_b(a)\,,$$
and this establishes our claim. \QED

It follows now readily from Lemma \ref{L5} that $\nabla^{(h)}_k$ and $\Sigma^{(h)}_k$ have the same asymptotic behavior. Specifically, we have:
\begin{corollary} \label{C1}
In the regime where $h\to \infty$ with $k^2\ll h \ll d^k$ and $\{\log_d h\}\to b\in[0,1)$, there is the weak convergence 
$$h d^{-k} \nabla^{(h)}_k - c(k-\lfloor \log_d h\rfloor)  \ \Longrightarrow \ L_b(c).$$
\end{corollary}
\proof Indeed we get  from Lemma \ref{L5} that for $h\ll d^k$
$$\E\left(\Sigma^{(h)}_k-\nabla^{(h)}_k\right) = O\left( d^k k^2 h^{-2}\right)\,,$$
and therefore if further $k^2\ll h$, then  
$$\lim h d^{-k} \left(\Sigma^{(h)}_k-\nabla^{(h)}_k\right)=0\qquad \hbox{in probability.}$$
We can thus conclude from Lemma \ref{L6}. \QED

We can now complete the proof of Theorem \ref{T2} along the same line as for Theorem \ref{T1};
for the sake of avoiding repetitions, technical details will be skipped.

\noindent{\bf Proof of Theorem \ref{T2}: \hskip10pt}    
For every $h\geq 1$, we write $(W^{(h)}_j, j\in \N) $ for a Galton-Watson process with binomial reproduction law with parameter $(d,\e^{-c/h})$. Then $\left(d^{-j}\e^{cj/h}W^{(h)}_j: j\in\N\right)$ is a square integrable martingale which converges a.s., and more precisely, one readily checks from martingale arguments that for every $\alpha >1/2$
$$\lim_{n\to \infty}Ê\P_n\left(\sup_{j\geq 0} \left|d^{-j}\e^{cj/h}W^{(h)}_j - n\right| > n^{\alpha}\right)=0
\qquad \hbox{uniformly in }h\geq 1\,,$$
where the notation $\P_n$ refers to the law of the process $W^{(h)}$ started from $W^{(h)}_0=n$ ancestors. 

We fix a starting level $k=\lfloor \log_d^4 h\rfloor$, and consider the process which counts for $j=0,1, \ldots$ the number of vertices of the $d$-regular tree at level $k+j$ which are still connected to the root after percolation. We obtain a version of the Galton-Watson process above, starting from $W^{(h)}_0= d^k-\nabla^{(h)}_k$. 
Skorokhod's representation theorem enables us to assume that the weak convergence in Corollary  \ref{C1} holds almost surely, that is 
$$h d^{-k} \nabla^{(h)}_k = c(k-\lfloor \log_d h\rfloor)  + L_b(c) + o(1).$$

We now observe the process for $j=h-k$, so $W^{(h)}_{h-k}= d^h-\nabla_h$,
and we deduce from above that  for e.g.  $\alpha=2/3$, we have with high probability
\begin{eqnarray*}
d^{k-h}\e^{c(h-k)/h}(d^h-\nabla_h) &=& d^k-\nabla^{(h)}_k + o(d^{2k/3})\\
&=& d^k-\frac{d^k}{h}c(k-\lfloor \log_d h\rfloor)  -\frac{d^k}{h} L_b(c) + o(d^k/h). 
\end{eqnarray*}
It follows that 
\begin{eqnarray*}
1-d^{-h}\nabla_h&=& \e^{-c(1-k/h)} \left(1- \frac{c}{h}(k-\lfloor \log_d h\rfloor)
-\frac{L_b(c)}{h} + o(1/h)\right) \\
&=& \e^{-c} + \frac{c\e^{-c} }{h}\lfloor \log_d h\rfloor 
-\e^{-c}\frac{L_b(c)}{h} + o(1/h) \,
\end{eqnarray*}
and this entails our claim. \QED

\end{section}

  \end{document}